\documentclass[12 pt]{article}
\usepackage{hyperref,amsmath,amsfonts,graphicx,psfrag,fancybox,qsymbols,indentfirst}
\usepackage[latin1]{inputenc}
\begin{document}
\setcounter{section}{0}
\title{Canonical moments and random spectral measures}
\large
\author{F. Gamboa  and A. Rouault}
\setcounter{tocdepth}{3}
\maketitle
\normalsize\rm
\begin{abstract}
We study some connections between 
the random moment problem and the random matrix theory. A uniform draw in a space of moments can be lifted into the spectral probability measure of the pair $(A, e)$ where $A$ is a random matrix from a classical ensemble and $e$ is a fixed unit vector. This random measure is a weighted sampling among the eigenvalues of $A$. We also study the large deviations properties of this random measure when the dimension of the matrix increases. 
The rate function for these large deviations  involves the reversed Kullback information.
\end{abstract}
\def \el{\sur{=}{(d)}}
\def \be{\begin{eqnarray*}}
\def \ee{\end{eqnarray*}}
\def \ben{\begin{eqnarray}}
\def \een{\end{eqnarray}}
\def\wh{\widehat}
\def\wt{\widetilde}
\def\pn{\par\noindent}
\def\sn{^{(N)}}
\def\AArm{\fam0 }
\def\AAk#1#2{\setbox\AAbo=\hbox{#2}\AAdi=\wd\AAbo\kern#1\AAdi\varphi}%
\def\AAr#1#2#3{\setbox\AAbo=\hbox{#2}\AAdi=\ht\AAbo\raise#1\AAdi\hbox{#3}
}%
\def \eref#1{(\ref{#1})}
\def \sur#1#2{\mathrel{\mathop{\kern 0pt#1}\limits^{#2}}}
\def\BBa{{\AArm A}}%
\def\BBb{{\AArm I\!B}}%
\def\BBd{{\AArm I\!D}}%
\def\BBe{{\AArm I\!E}}%
\def\BBf{{\AArm I\!F}}%
\def\BBg{{\AArm G\AAk{-1.02}{G}\AAr{.9}{I}{\AAFf\char"3F}}}%
\def\BBh{{\AArm I\!H}}%
\def\BBi{{\AArm I\!I}}%
\def\BBj{{\AArm J}}%
\def\BBk{{\AArm I\!K}}%
\def\BBl{{\AArm I\!L}}%
\def\BBm{{\AArm I\!M}}%
\def\BBn{{\AArm I\!N}}%
\def\N{{\AArm I\!N}}%
\def\BBo{{\AArm O\AAk{-1.02}{O}\AAr{.9}{I}{\AAFf\char"3F}}}%
\def\BBp{{\AArm I\!P}}%
\def \dis{\displaystyle}
\def\P{{\AArm I\!P}}%
\def \p{{\cal P}}
\def\BBq{{\AArm Q\AAk{-1.02}{Q}\AAr{.9}{I}{\AAFf\char"3F}}}%
\def\Q{{\AArm Q\AAk{-1.02}{Q}\AAr{.9}{I}{\AAFf\char"3F}}}%
\def\BBr{{\AArm I\!R}}%
\def\BBs{{\AArm S}}%
\def\BBu{{\AArm U\AAk{-1}{U}\AAr{.95}{I}{\AAFf\char"3F}}}%
\def\BBv{{\AArm V}}%
\def\BBw{{\AArm W}}%
\def\BBx{{\AArm X}}%
\def\BBy{{\AArm Y}}%
\def\BBz{{\AArm Z\!\!Z}}%
\def\BBone{{\AArm 1\AAk{-.8}{I}I}}%
\def\D{{\cal D}}
\def\R{{\math R}}
\def\N{{\math N}}
\def\E{{\math E}}
\def\P{{\math P}}
\def\Z{{\math Z}}
\def\Q{{\math Q}}
\def\C{{\math C}}
\def \e{{\rm e}}
\def \f{{\cal F}}
\def \g{{\cal G}}
\def \L{{\cal L}}
\def \d{{\tt d}}
\def \k{{\tt k}}
\def \i{{\tt i}}
\def \w{{\tt w}}
\def \s{{`S`O}_{2N}}
\def\rp{\right)}
\def\lp{\left(}
\def\Brp{\Big)}
\def\Blp{\Big(}
\def\brp{\big)}
\def\blp{\big(}
\def\noi{\noindent\underbar}
\def \imp{\Rightarrow}
\def\cb{\hbox{C}\beta\hbox{E}(N)}
\def\so{`S`O(2N)}
\def\E{\ensuremath{\mathbb{E}}}
\newcommand{\ed}{\mbox{$ \ \stackrel{d}{=}$ }}
\def\videbox{\mathbin{\vbox{\hrule\hbox{\vrule height1ex \kern.5em\vrule height1ex}\hrule}}}
\font\calcal=cmsy10 scaled\magstep1
\def \build#1#2#3{\mathrel{\mathop{\kern 0pt#1}\limits_{#2}^{#3}}}
\def\suPP#1#2{{\displaystyle\sup _{\scriptstyle #1\atop \scriptstyle #2}}}
\def\proDD#1#2{{\displaystyle\prod _{\scriptstyle #1\atop \scriptstyle #2}}}
\def\proof{\noindent{\bf Proof:}\hskip10pt}
\def\QED{\hfill\vrule height 1.5ex width 1.4ex depth -.1ex \vskip20pt}
\newcommand{\ii}{{\mathrm{i}}}
\newcommand{\cCn}{\mathcal{C}_n}
\newcommand{\br}{\mathcal{B}}
\newcommand{\usn}{a_{n}}
\newcommand{\Int}{\operatorname{int}\,}
\newcommand{\Clo}{\operatorname{clo}\,}
\newcommand{\prs}[2]{\langle#1,#2\rangle}
\newcommand{\eul}{\operatorname{e}}
\newcommand{\dd}{\operatorname{d}}
\newcommand{\tr}{\operatorname{tr}}
\def\cro#1{\llbracket #1  \rrbracket}

\newtheorem{lem}{Lemma}[section]
\newtheorem{defi}[lem]{Definition}
\newtheorem{theo}[lem]{Theorem}
\newtheorem{cor}[lem]{Corollary}
\newtheorem{prop}[lem]{Proposition}
\newtheorem{nota}[lem]{Notation}
\newtheorem{rem}[lem]{Remark}
\newtheorem{rems}[lem]{Remarks}
\newqsymbol{`P}{\mathbb{P}}
\newqsymbol{`E}{\mathbb{E}}
\newqsymbol{`B}{\mathbb{B}}
\newqsymbol{`N}{\mathbb{N}}
\newqsymbol{`C}{\mathbb{C}}
\newqsymbol{`D}{\mathbb{D}}
\newqsymbol{`R}{\mathbb{R}}
\newqsymbol{`S}{\mathbb{S}}
\newqsymbol{`O}{\mathbb{O}}
\newqsymbol{`o}{\omega}
\newqsymbol{`M}{{\cal M}}
\newqsymbol{`H}{\mathbb{H}}
\newqsymbol{`X}{\mathbb{X}}
\newqsymbol{`T}{\mathbb{T}}
\newqsymbol{`U}{\mathbb{U}}
\newcommand{\CUE}{\operatorname{CUE}}
\newcommand{\GUE}{\operatorname{GUE}}
\newcommand{\GOE}{\operatorname{GOE}}
\newcommand{\CSE}{\operatorname{CSE}}
\newcommand{\COE}{\operatorname{COE}}
\newcommand{\Beta}{\operatorname{Beta}}
\newcommand{\ARCSIN}{\operatorname{arcsine}}
\newcommand{\USp}{\operatorname{USp}}
\newcommand{\Dir}{\operatorname{Dir}}
\newcommand{\BBt}{\mathbb{T}}
\newcommand{\BBc}{\mathbb{C}}
\newcommand{\diag}{\operatorname{diag}}
\newcommand{\CbetN}{\operatorname{C\beta E}_{N}}
%
%
\newdimen\AAdi%
\newbox\AAbo%
\font\AAFf=cmex10

\setcounter{page}{1}

\def \build#1#2#3{\mathrel{\mathop{\kern 0pt#1}\limits_{#2}^{#3}}}
\def\suPP#1#2{{\displaystyle\sup _{\scriptstyle #1\atop \scriptstyle #2}}}
\def\proDD#1#2{{\displaystyle\prod _{\scriptstyle #1\atop \scriptstyle #2}}}
\def\proof{\noindent{\bf Proof:}\hskip10pt}
\setcounter{page}{1}
\setcounter{section}{0}

\section{Introduction}
\label{suintr}

The aim of this paper is to establish a strong link between
random moment problems (\cite{Chang},\cite{FabLoz}) and random matrices.

As a matter of fact, in the last decade, emphasis has been put on the
asymptotic behavior of large random matrices. Such a study  
was first motivated by problems arising in theoretical physics (\cite{Mehta}, \cite{Forrester}).
Generally, the distribution of these random matrices is characterized  by some invariance properties. In  the $\GOE_{N}$ (resp. $\GUE_{N}$),  the matrix $A_N$ is $N \times N$ symmetric (resp. Hermitian) and, except for this constraint, has independent Gaussian entries. In the $\CUE_{N}$, $A_N$ has the Haar distribution on $`U(N)$ (the group of $N\times N$ unitary matrices). We refer to
\cite{Baimethodo} and \cite{HiaiP} for a general overview on the subject. In these examples, the distribution of both the random eigenvalues $\lambda_1,\ldots,\lambda_N$ and eigenvectors are precisely known.  
  These eigenvalues are (stochastically) independent of the eigenvectors and the matrix of the normalized eigenvectors is Haar distributed. 
%
The main object of interest is the empirical  spectral distribution (ESD) of $A_N$,
which is the discrete measure
\ben
\nu\sn:=\frac{1}{N} \sum_{j=1}^{N}\delta_{\lambda_j}. 
\label{esd}
\een
When  
 $A_N$ is  drawn randomly,   
  the asymptotic properties of the sequence of random probability measures $(\nu^{(N)})_N$ 
are obtained by a moment method using 
\begin{equation}
\int x^{k}\nu\sn(\dd\! x)=\frac{1}{N}\tr(A_N^k),\; (k\geq 1)\,,
\label{tratra}
\end{equation}
%
%
 and 
$(\nu^{(N)})_N$ converges towards a deterministic probability measure. Namely, for the $\GOE$ (resp. $\GUE$) the limit is the semicircle distribution while for the $\CUE$, it is the uniform measure on the unit circle $\mathbb T$. 
To give 
 the rate of convergence of the sequence $(\nu\sn)_N$, 
it is interesting to know
 if it satisfies the large deviation principle (LDP).
This problem has been first studied for the $\GOE$ in the pioneer work of Ben Arous and Guionnet (\cite{BenGui}). It was further investigated 
for the CUE (see \cite{HiaiIHP})  and many other distributions (see  
\cite{EichSto}).
In these LDPs 
 the speed is $N^{2}$ and 
  the rate function may be obtained from the Voiculescu entropy (\cite{HiaiP}).

In this paper, we consider another approach to unitary matrices. 
Instead of considering the classical empirical measure $\nu^{(N)}$ defined in (\ref{esd}), we  study the 
 random spectral measures (see (\ref{wsd}) below).
Every unitary operator $A$ on a Hilbert space with a unit cyclic vector $e$ is (by the spectral theorem) unitarily equivalent to the multiplication by $z$ on $L^2(\mathbb T, \mu)$ for some probability measure $\mu$, which is uniquely determined by its moments.
If we  endow $\mathbb U (N)$  with the normalized Haar measure $\lambda_{\mathbb U (N)}$, the first vector of the standard basis $e_1$ is a.s. cyclic for $A_N$. We thus define  the random spectral  measure $\mu\sn_{\w}$ on $\mathbb T$ related to the pair $(A_N, e_1)$. 
It is defined by its moments:  
\begin{equation}
\int_{\mathbb T} z^{k}\mu\sn_{\w}(\dd\! z)=\langle e_1, A_N^k e_1\rangle,\; (k\geq 1)\,.
\label{tratra1}
\end{equation}
The ESD $\nu^{(N)}$ and the spectral measure $ \mu\sn_{\w}$ are widely different. This can be seen in terms of moments.
In both cases,  we consider successive powers of $A_N$, but the moments of $ \mu\sn_{\w}$ are the first entries given by (\ref{tratra1}) and those of $\nu^{(N)}$  are the normalized traces given by 
(\ref{tratra}). 
Alternatively, a simple diagonalization of $A_N$ gives:  
\ben
\mu_\w\sn:=\sum_{j=1}^{N}\w_j\delta_{\eul^{i\theta_j}}\,, 
\label{wsd}
\een
where the $\w_j$'s are the square moduli of the top entries of the normalized unit eigenvectors and 
$\exp(i\theta_j)=\lambda_j$ are the (a.s. different) eigenvalues of $A_N$. 
It is now clear that
 this measure carries more information on $A_N$ than $\nu^{(N)}$. Ideed,  it involves weights depending on the eigenvectors of this matrix. 
%
%
In this paper, we are concerned with a precise study of the random probability measure $\mu_\w\sn$.
%

For a probability measure $\xi$ on $\mathbb T$, two 
encodings  are interesting.  First, its Verblunsky coefficients appear in the Schur recursion formula satisfied by the sequence of orthogonal polynomials in $L^2(\xi)$ (see for example \cite{Simon1}, \cite{Simon2}). Secondly, its canonical moments are involved in statistical applications.
They are inductively defined : the $(k+1)$-th canonical 
moment is the relative position of the $(k+1)$-th moment 
consistent with the $k$ first ones
(see \cite{DS} Chapter 9 for an 
overview on canonical moments).

Two remarkable papers motivated our work. On the one hand, 
Killip and Nenciu \cite{Killip1} 
proved 
 that in the $\CUE_N$ model, the Verblunsky coefficients of 
$\mu_\w\sn$ are independent and have known distributions (namely complex beta distribution). 
On the other hand, Lozada \cite{LozEJP} has studied the distribution of the $N$ first canonical moments 
for special random measures, namely the measures having
their $N$ first moments uniformly distributed.
Surprisingly, the distributions found by Killip and Nenciu and by Lozada are the same. We can explain this last result  by the fact that Verblunsky coefficients and canonical moments coincide
(as pointed out by several authors, see Simon \cite{Simon3} p. 439). Using this observation we show in Theorem \ref{=C} that the vector of moments of
$\mu_\w\sn$ under $\CUE_{N}$  is also uniformly distributed, giving a (new) enlightening connection between random moment theory and random matrix theory. 

%
%
%
%
Extending this method, we give 
 the precise distributions of moments and canonical moments 
 in other circular ensembles and in  the so-called $\beta$-ensembles (or log-gases). 
Furthermore, 
we consider real matrices. Beginning with  the special orthogonal group of order $2N$ $`S`O(2N)$, the symmetry property of spectral measures allows to consider their projections on $[-1 , 1]$, which leads 
 naturally to the family of Jacobi ensembles. We prove a lift of a uniform random point in the space of moments into the spectral measure of a random element of the ensemble $\mathbb S\mathbb O (2N)/\mathbb U(N)$ (Theorem \ref{unif2}).

The common feature of all these models 
 is that the weights and the support points of $\mu_\w\sn$ are independent. It is a nice framework 
to  study the LDP, 
with  two main differences with the above mentioned  LDP for ESD : the (random) weighting  both slows down the speed from $N^2$ to $N$ and (consequently) changes the rate function to the reversed Kullback information. 
This representation
allows in particular to recover  results proved by using the independence of canonical moments in \cite{FabLoz} and \cite{LozEJP}.
As a consequence, we see that $\mu_\w\sn$ converges weakly in probability to the same deterministic limits as the ESD $\nu\sn$.

Some results for ensembles with non compactly supported spectral measures  can be found in \cite{Bai07} for sample covariance matrices, and in \cite{gr+} for the GUE.

The paper is organized as follows.
In Section \ref{sunodef}, we give some notations and list some  distributions 
 of frequent use. In Section \ref{cadixxx}, we explain the set-up of moments and canonical moments in the complex and real case. Section \ref{fontenay} is devoted to the distribution of random spectral measures in different models. Finally, in Section \ref{arld}, we state LDPs 
  for the families of random spectral measures.

Notice that  a few months after dropping the first version of the present
paper on Arxiv, we have been aware of the work of Birke and Dette \cite{birke}
on the asymptotic behavior of the roots of random orthogonal
polynomials associated to random moments. 
Their main result 
 is the computation of the root distribution 
when the 
 moments are uniformly distributed. 
 It appears
to be a particular case of our more general Lemma \ref{oneof} and Theorem \ref{bizth}.

\section{Notations and some useful distributions}
\label{sunodef}
Let $`D$ be the open unit disk in $`C$ and  $`T=\partial `D$, the unit circle parameterized by $z = \eul^{i\theta}$, with $\theta \in [-\pi, \pi)$. Let $\mu$ be a probability measure on $`T$.  $\mu$ is said to be trivial if it supported on a finite set, and nontrivial otherwise. Further, the Hermitian product on $L^2(`T,\mu)$ is defined by 
$\langle f , g\rangle:=\int_{`T} \bar{f}\; g\;\dd\! \mu$.
Let $(G,\mathcal{G})$  be any measurable space, we denote by ${\mathcal M}_1(G)$ 
(resp. ${\mathcal M}(G)$) the set of all probability measures (resp. positive measures) on $G$.
If $\mu\in{\mathcal M}_1(G)$ and $f$ is integrable, we write sometimes $\mu(f)$ for 
$\int_G f \dd\! \mu$. 
We  recall now some special useful distributions that are
used later. 
For $k \geq  1$, we set
\be {\mathcal S}_{k} := \{(x_1, \cdots, x_{k}) : x_i > 0 , (i = 1, \cdots , k) ,  \ x_1 + \cdots + x_{k}  = 1 \} \\
{\mathcal S}_k^< := \{(x_1, \cdots, x_k) : x_i > 0 , (i = 1, \cdots , k) ,  \ x_1 + \cdots + x_k  < 1 \}\,.
\ee
Obviously, the mapping $(x_1, \cdots , x_{k+1}) \mapsto (x_1, \cdots, x_k)$ is one to one from the simplex ${\mathcal S}_{k+1}$ onto  ${\mathcal S}_k^<$.

For $a_j>0,\;j=1,\ldots, k+1$,
the Dirichlet distribution $\Dir(a_1, \cdots, a_{k+1})$ on ${\mathcal S}_{k+1}$ 
has the density
\ben
\label{defdir} 
\frac{\Gamma(a_1 + \cdots + a_{k+1})}{\Gamma(a_1) \cdots \Gamma(a_{k+1})}\!\ x_1^{a_1 - 1} \cdots x_{k+1}^{a_{k+1} -1}
\een
with respect to the Lebesgue measure on ${\mathcal S}_{k+1}$. When $a_1= \cdots = a_{k+1} = a>0$, we denote the Dirichlet distribution by
$\Dir_{k+1} (a)$. If $a=1$ we recover  the uniform distribution on ${\mathcal S}_{k+1}$.

Pushing the Dirichlet distribution under the previous mapping, we get a probability measure  $\Dir(a_1, \dots,a_k;a_{k+1})$ on 
  ${\mathcal S}_k^<$ with density :
\be
\frac{\Gamma(a_1 + \cdots + a_{k+1})}{\Gamma(a_1) \cdots \Gamma(a_{k+1})}\!\  x_1^{a_1 - 1} \cdots x_k^{a_k -1} \left(1 - x_1 - \cdots -x_k\right)^{a_{k+1}-1}
\label{defdir2}
\ee

If $k=1$,  $\Dir(a ; b)$ has the  $\Beta(a,b)$  density on $[0,1]$ :
$$\frac{\Gamma(a+b)}{\Gamma(a)\Gamma(b)} x^{a-1}(1-x)^{b-1} 
\,.$$
The very particular case $\Beta(1/2,1/2)$ is the so-called arcsine law. 
Sometimes we need the distribution obtained by pushing forward $\Beta(a,b)$ under the mapping $x \mapsto 2x-1$. 
It is the distribution $\Beta_s(b,a)$ on $(-1,1)$ having density
$$2^{1-a-b}\frac{\Gamma(a+b)}{\Gamma(a)\Gamma(b)} (1+x)^{a-1}(1-x)^{b-1}
\,.$$
We use also a complex version of the beta distribution. For $r > -1$ let $\eta_r$ be the probability density on $\mathbb C$ defined by
\ben\label{eta}
\eta_r (z):= \frac{r+1}{\pi}\!\ \big(1- |z|^2\big)^r\;,\; (z\in \mathbb D).
\een
It is obviously the density of $X= \eul^{\ii U}\sqrt B$ where $U$ is uniform on $[0,2\pi]$ and $B$ is $\Beta(1, r)$ distributed.

To end this section, let us recall the classical relation between the Dirichlet and Gamma distributions.  
A Gamma variable $\gamma_a$ with parameter $a > 0$ has density
\[\frac{t^{a-1} \eul^{-t}}{\Gamma(a)}, \ (t> 0). \]
It is well known that if $y_i, i =1, \cdots , r$ are independent and if $y_i \el \gamma_{b_i}$ ($b_i>0$) then
\ben
\label{10}
\left(\frac{y_1}{y_1+\cdots +y_{r}} , \cdots , \frac{y_{r}}{y_1+\cdots +y_{r}}\right) \el \Dir\left(b_1 , \cdots , b_{r}\right)
\een
and this variable is independent of  $y_1+\cdots +y_{r}$.
\section{Verblunsky coefficients, canonical moments and uniform probability}
\label{cadixxx}
We now present the notion of canonical moments of a measure. We successively focus on the complex case (unitary circle) and on the real case (compact interval). Finally we recall that the uniform probability on moment spaces corresponds to a nice distribution on the space of canonical moments. 
\subsection{The complex case : Verblunsky coefficients}
\label{vertsky}
All the material of this subsection comes from \cite{Simon1} Section 1 or \cite{Simon3}
Sections 2 and 3. We recall here the connection between moments of a probability
measure on $`T$ and Verblunsky coefficients built through orthogonal polynomials.

Let $\mu$ be an arbitrary 
nontrivial probability measure on $`T$. 
The functions $1, z, z^2, \cdots$ are linearly independent in $L^2(`T, \mu)$. 
Following the Gram-Schmidt procedure we define the infinite sequence  $(\Phi_j)_{j \geq 1}$ of monic
orthogonal polynomials. More precisely, $\Phi_0 (z) \equiv 1$ and 
$\Phi_j (z)$ 
is the projection of $z^j$ onto $\{1, \cdots,z^{j-1}\}^\perp$, for $j\geq 1$. 

If $\mu\in {\mathcal M}_1(`T)$ has finite support $\{z_1, \cdots , z_K\}$ ($K$ different points), we still define $\Phi_j$ in the same way for $j= 1, \dots, K-1$. Besides, 
we define $\Phi_K$ as the unique monic polynomial of degree $K$ such that 
$\Vert \Phi_K\Vert =0$ i.e.
\ben\label{polycar}\Phi_K (z) = \prod_{j=1}^K (z-z_j)\,.\een
By convention, in the nontrivial case we set $K=\infty$, and then the sentence "for every $j < K$" (or "for every $j < K+1$") mean "for every finite $j$".
 
Some useful polynomials associated to the sequence $(\Phi_j)_{j < K}$ are
the reversed (or reciprocal) polynomials. They are defined by $\Phi_0^\star (z)\equiv 1$ and 
\ben
\label{reversed}
\Phi_j^\star (z) = z^j \overline{\Phi_j(1/\bar z)}\,.
\een
In the next proposition, we define the Verblunsky coefficients. They are sometimes called 
Schur, Szeg\H{o} or Geronimus coefficients or even reflection coefficients. 
Surprisingly, these coefficients also appears as  central quantities in the theory of moment problems
(see for example \cite{DS}), where they are called canonical moments.
\begin{prop}[Szeg\H o]
\label{propcan}
For $j< K+1$, we define the Verblunsky coefficient of order $j$ by setting
$c_j:=-\overline{\Phi_{j}(0)}$. 
The sequence $(\Phi_j)_{j < K+1}$
  satisfies the  recursion
\ben
\label{recursion}
\nonumber
\Phi_0 = 1 \ \ , \ \ \Phi_{j} = z \Phi_{j-1} - \bar c_j \Phi_{j-1}^\star\ \ \ (1 \leq j < K+1)\,.
\een
\end{prop}
Notice that if $\mu$ is nontrivial, $c_j\in \mathbb D$ for every $j> 0$. 
If  $K < \infty$, then $c_j\in \mathbb D$ for $1\leq j\leq K-1$ and $c_{K} \in `T$  
(see \cite{Simon1} Theorem 1.5.2). In the sequel, when we need
to stress the dependence
of the Verblunsky coefficients  on the underlying measure $\mu$, 
we  
  write $c_j(\mu)$. 
%
\label{211}
A theorem due to Verblunsky 
claims that the  correspondence between  a probability measure 
 on $\mathbb T$ and the sequence of its coefficients is one to one (\cite{Simon1} Theorem 1.7.11).
%
For $N \geq 1$, set 
\ben
M_N^{`T} = \left\{\left(\int_{`T} z^j \mu(\dd\! z)\right)_{1\leq j\leq N}:\; \mu \in {\mathcal M}_1(`T)\right\}\,.
\een
\begin{prop}[\cite{Simon1} pp. 60 and 218, \cite{DS} p.269]
Let $\left(t_1, \cdots , t_N\right)\in \Int M_N^{`T}$, the range of the $(N+1)$th moment 
$$t_{N+1} = \int_{`T} z^{N+1}  \eta(\dd\! z)$$
as $\eta$ varies over all probability measures having $(t_1, \cdots , t_N)$ as their $N$ first moments,
is a disk centered at some point $s_{N}$ (depending on $t_1, \cdots , t_N$)  
with radius 
$$r_{N} = \prod_{j=1}^{N}\left(1 -|c_j|^2\right)\not= 0\,$$
The relative position is
$$\frac{t_{N+1} - s_{N}}{r_{N}}\in `D\,,$$
and it
is exactly $\bar c_{N+1}$.
\end{prop}
This striking connection of two notions (related position versus recursion coefficient) has been proved  in two opposite ways: starting from the recursion to get relative position \cite{Simon1} p. 60), or alternatively starting from the moments and their relative positions and computing the recursion  \cite{Simon1} p. 218,  \cite{DS} p.269. 
For an historical account, see \cite{Simon1} p.10 and p.221.

%
\subsection{Canonical moments : Real case}
\label{2.2}
Following the last interpretation of canonical moments as relative positions of moments in
$M_N^{`T}$,
it is possible to build similar quantities when the set of integration $`T$ is replaced by a compact interval $[a,b]$ ($a<b$). Let us briefly recall the exact definition of canonical moments in this frame. We refer to the excellent book of Dette and Studden \cite{DS} for a complete overview on the subject. 
To begin, 
for $N\geq 1$,  we denote  by $M_N^{[a,b]}$ the $N$-th moment space generated by probability measures on $[a,b]$:
\[M_N^{[a,b]} := \left\{\left(\int_a^b x^k  \mu( \dd\! x) , k =1, \dots, N\right)\ : \; \mu \in {\mathcal M}_1([a,b])\right\}\,.\]

Given ${\bf m}^j := (m_1, \cdots ,m_j) \in \Int M_j^{[a,b]},\;
(j \geq 1)$, we first define the {\it extreme} values, 
\ben
m_{j+1}^+({\bf m}^j) &=& \hbox{max}\left\{m \in `R : (m_1, \cdots , m_j, m) \in M_{j+1}^{[a,b]}\right\}
\label{doudou}
\\
m_{j+1}^-({\bf m}^j) &=& \hbox{min}\left\{m \in `R : (m_1, \cdots , m_j, m) \in M_{j+1}^{[a,b]}\right\}.
\label{doudou1}
\een
For $k\geq i \geq 1$, the $i$-th canonical moment is defined recursively as
\ben
\label{defcanonr}
c_i ({\bf m}^k) = c_i ({\bf m}^i) := \frac{m_i - m_{i}^-({\bf m}^{i-1})}{ m_{i}^+({\bf m}^{i-1})- m_{i}^-({\bf m}^{i-1})}
\een
A quite nice property of canonical moments 
is that they are invariant
on any affine one to one mapping transforming the support of the underlying measures (see for example \cite{DS}).
So that,  we  may restrict ourselves to the special case of $M_N^{[0,1]}$.

\subsection{Uniform probability on moment spaces}
In this subsection, we  endow the moment sets $M_n^{`T}$ and $M_n^{[0,1]}$ (for fixed $n$) with the uniform probability and recall that in these cases the canonical moments previously defined have very interesting properties (see Lemmas \ref{unifreal} and \ref{unifcomplex} below).

As pointed in Section \ref{vertsky},  for fixed $N$ there is a one to one mapping between moments and canonicals moments. 
Indeed, we may define a one to one mapping $\kappa_N^{\mathbb T}$
\ben\nonumber \kappa_N^{\mathbb T} \ : \ \ \Int M_N^{\mathbb T} &\rightarrow& {\mathbb D}^N\\
(t_1, \dots, t_N) &\mapsto& (c_1, \dots, c_N)\,.
\een
In the real case, the relation (\ref{defcanonr}) defines a one to one mapping $\kappa_N^{[0,1]}$
\ben\nonumber \kappa_N^{[0,1]} \ : \ \ \Int M_N^{[0,1]} &\rightarrow& (0,1)^N\\
(t_1, \dots, t_N) &\mapsto& (c_1, \dots, c_N)\,.
\een
which is {\it triangular} and bijective.


The two following lemmas give the canonical moment distribution when the moments are uniformly drawn.

\begin{lem}[Lozada]
\label{unifcomplex}
Endowing $\Int M_N^{`T}$ with the uniform distribution is equivalent to the $N$ first canonical moments $(c_1, \dots, c_N)$ being independent in such a way that $c_j$ has density $\eta_{N-j}$.
\end{lem}

\begin{lem}[Chang, Kemperman, Studden]
\label{unifreal}
Endowing $\Int M_N^{[0,1]}$ with the uniform distribution is equivalent to the $N$  canonical moments $(c_1, \dots, c_N)$ being independent in such a way that $c_j$  is $\Beta(N-j+1, N-j+1)$ distributed.
\end{lem}
These two lemmas have been the starting point for the investigation on the asymptotic behavior of the randomized sets of moments $M_N^{[0,1]}$ and $M_N^{`T}$ (see
\cite{Chang}, \cite{FabLoz}, \cite{LozEJP} and \cite{gamb}). Here, these lemmas are useful to obtain some very nice elementary properties of random measures built on eigenvalues of some classical random matrix models. We develop these results in the next sections.  
\section{Random spectral measures and their moments}
\label{fontenay}
In this Section, we first define the spectral measure associated with the pair $(A,e)$ where $A$ is a unitary matrix and $e$ a cyclic vector. Then 
we randomly draw a matrix and study the distribution of the associated random spectral measure and of its canonical moments. 
We  focus on various classical distributions popular in the random matrix paradigm. Moreover we extend the class of such random measures leading to log-gases models, which in turn may be lifted into matrix ensembles. In particular, we emphasize the different ways to get uniform distribution on the space of moments. Our main results in this section are Theorem \ref{=C}, Theorem \ref{bizth} and Theorem \ref{unif2}.


\subsection{Spectral measures associated with a unitary matrix}

Let us consider ${\mathbb C}^N$ with its canonical basis $(e_1, \dots, e_N)$ and ${\mathbb U}(N)$ the group of $N\times N$ unitary matrices. For $A \in `U(N)$, let $\lambda_i,\; i=1,\ldots,N$ be the eigenvalues of $A$.  We may write 
\begin{equation}
A=\Pi D\Pi^*, 
\label{march}
\end{equation}
where $D:=\diag(\lambda_i)_{i=1,\ldots,N}$
is diagonal  and $\Pi:=(\pi_{ij})_{i,j=1,\ldots,N}$ is unitary. Let us assume that $e_1$ is cyclic, i.e. that span$(e_1, Ae_1, \dots, A^{N-1}e_1)$ has rank $N$. Following \cite{Simon2}, we consider the spectral
measure $\mu_\w\sn$ associated with the pair $(A, e_1)$ i.e. having, for any $k\in\BBn$, moment of order $k$ equal to $\prs{e_1}{A^ke_1}$.
Obviously, $\mu_\w\sn$ is trivial and supported by the eigenvalues of $A$. More precisely,
 \ben
\label{mastermu}
\mu_\w\sn = \sum_k \w_k
 \ \! \delta_{\lambda_k
}\,,\een  
where $\w_k := |\pi_{1k}|^2, k=1, \dots, N$ 
are  the square moduli of the top entries of the normalized  eigenvectors. 
\subsection{Circular ensembles}
Here are recalled three families of classical ensembles of random matrices:
\begin{itemize}
\item $\CUE_{N}$ (circular unitary ensemble) :  $\mathbb U(N)$ equipped with its normalized Haar measure $\lambda_{\mathbb U (N)}$
\item $\COE_{N}$ (circular orthogonal ensemble) : set of symmetric unitary $N\times N$ matrices. Every element $S$ can be written $S = g^T g$ with $g\in \mathbb U(N)$ hence $\COE_{N} \cong \mathbb U(N)/\mathbb O(N)$. The probability measure is the pushforward of $\lambda_{\mathbb U (N)}$  under the projection mapping $\mathbb U(N) \rightarrow \mathbb U(N)/\mathbb O(N)$.
\item $\CSE_{N}$ (circular symplectic ensemble) : set of $2N \times 2N$ self-dual unitary matrices. Recall that the dual of a matrix $H$ is $H^D := JH^T J^T$ with 
\[J = \begin{pmatrix} 0&-I_N\\ I_N&0\end{pmatrix}\,\]
and that a  $2N \times 2N$ matrix $k$ is symplectic if it satisfies $kJk^T = J$. 
Every element $\widetilde S$ can be written $\widetilde S = g^D g$ with $g\in \mathbb U(2N)$,
 hence $\CSE_{N} \cong \mathbb U(2N)/\USp(2N)$, where $\USp(2N)$, the set of unitary symplectic matrices, is the invariant set of the involution $g \mapsto (g^D)^{-1}$. The probability measure is the pushforward of $\lambda_{\mathbb U (2N)}$ under the projection mapping $\mathbb U(2N) \rightarrow \mathbb U(2N)/ \USp(2N)$.
\end{itemize}
All these ensembles are considered in the two following subsections.
\subsubsection{The circular unitary ensemble}
Let us first notice that  with probability $1$, the vector $e_1$ is cyclic.
If $A \in `U(N)$,  all its eigenvalues  have unit modulus and in the decomposition (\ref{march}) we may write $D = \diag \left(\eul^{i\theta_1}, \cdots, \eul^{i\theta_N}\right)$. Further, the matrix $\Pi$ in (\ref{march}) is not uniquely determined. 
But, it is well known that it can be chosen such that, if $A$ is  Haar distributed then $\Pi$ is also Haar distributed and independent of $D$. In this case, the joint law of eigenvalues arguments (on $[0,2\pi]^N$) is 
\ben
\label{loitheta}\CUE_{N} (\dd\!\theta_1, \cdots, \dd\! \theta_N) := \frac{1}{N! (2\pi)^N}\ |\Delta(\eul^{i\theta_1}, \cdots , \eul^{i\theta_N})|^2\ \dd\! \theta_1 \cdots \dd\!\theta_N\ \,
\een
where $\Delta$ is the Vandermonde determinant
\ben
\label{defdelta}
\Delta(x_1, \cdots, x_N) = \prod_{1\leq j< k\leq N} \left(x_j - x_k\right)\,.
\een
Furthermore, since $\Pi$ is Haar distributed, its first row is uniformly distributed on the Euclidean sphere of $\BBc^N$. 
So that, setting $\w_k := |\pi_{1,k}|^2, (k=1,\ldots,N)$ the  
vector $(\w_1, \cdots , \w_N)$ is uniformly distributed on ${\mathcal S}_N$ (see for example Proposition 3.1 of \cite{Killip1}) i.e. 
\ben
\label{loimuk}
(\w_1, \dots, \w_N) \el \Dir_N(1)
\een
(see (\ref{defdir}) for the definition).

So, from (\ref{mastermu}) we have built the random probability measure on $\mathbb T$:
\ben
\label{mastermurandom}
\mu_\w\sn = \sum_{k=1}^N \w_k
 \ \! \delta_{\eul^{i\theta_k}
}\,,\een
supported by the eigenvalues of the matrix $A$ with an independent system of weights. \\
Up to our knowledge, the following result  was not known previously.
\begin{theo}
\label{=C}
Under $\CUE_{N}$, the distribution of 
$\left(\int_{`T} z^j \mu_\w\sn(\dd\! z)\right)_{1\leq j\leq N-1}$
is uniform on $M^{`T}_{N-1}$. 
\end{theo}
\proof
Let  $(c_1(\mu_\w\sn), \cdots, c_{N}(\mu_\w\sn))$ be the vector of $N$ first canonical moments of $\mu_\w\sn$.  
In Proposition 3.3 of \cite{Killip1} (see also
\cite{simon2006cmf} Section 11) it is proved  that the components of this vector are independent. Furthermore,
\ben
\label{cuecj}
\CUE_{N} (c_j(\mu_\w\sn) \in \dd\! z) = \eta_{N-j-1}(z) \dd\! z\,, \ j= 1, \cdots , N-1\,,
\een 
and $c_{N}$ is uniform on $\mathbb T$ (since the support of $\mu_\w\sn$ is a finite set, the last canonical moment belongs to $\mathbb T$). 
Considering now the pushforward of these distributions under $(\kappa_{N-1}^{\mathbb T})^{-1}$ the measurable one to one mapping from 
$`D^{N-1}$ to $\Int M_{N-1}^{`T}$ 
which map  canonical moments to moments we may conclude using Lemma \ref{unifcomplex}.\QED
\noindent
For the special unitary group $\mathbb S\mathbb U (N)$ we also have.
\begin{cor}
\label{petitcor}
Under the normalized Haar measure $\lambda_{\mathbb S\mathbb U(N)}$ on $\mathbb S\mathbb U (N)$,  
the distribution of  $\left(\int_{`T} z^j \mu_\w\sn(\dd\! z)\right)_{1\leq j\leq N-1}$
is uniform on $M^{`T}_{N-1}$.
\end{cor}
\proof 
The joint eigenvalue distribution on $\mathbb T^{N-1}$ of the normalized Haar measure on $\mathbb S\mathbb U (N)$ is
\[\frac{1}{N! (2\pi)^{N-1}}\ |\Delta(\eul^{i\theta_1}, \cdots , \eul^{i\theta_N})|^2\ \dd\!\theta_1 \cdots \dd\!\theta_{N-1}\ \,\]
with $\theta_N = -(\theta_1+ \dots + \theta_{N-1}) (\hbox{mod} \!\ 2\pi)$ 
(see \cite{hiai3}).
Notice that it is the  density of $(\theta_1, \dots, \theta_{N-1})$ under $\CUE_{N}$ conditioned on $\theta_N = -(\theta_1+ \dots + \theta_{N-1}) (\hbox{mod} \!\ 2\pi)$.
Moreover we have  
$\exp \big(i (\theta_1+ \dots + \theta_N)\big) = (-1)^N \Phi_N(0) = (-1)^{N+1} \bar c_N$,
so that the conditioning set is $\{c_N = (-1)^{N+1}\}$. We have seen in the proof of  Theorem
\ref{=C} that $c_1,\cdots, c_N$ are independent under $\CUE_{N}$. Therefore, $c_1,\cdots, c_{N-1}$ remain independent when conditioning on  
$\{c_N = (-1)^{N+1}\}$. Moreover, their distributions are not affected by this conditioning. 
 Using Lemma \ref{unifcomplex} we deduce that the conclusion of the Corollary holds true.
\QED
\subsubsection{The circular $\beta$-models}
\label{goupil}
For the $\COE$ and $\CSE$ 
 we have 
 also independence between the weights and the support points. Moreover, the distributions are given by  formulas analogous  to (\ref{loitheta}) and (\ref{loimuk}).  Actually it is convenient (and now classical) to define a class of distributions with a continuous parameter $\beta > 0$, extending the cases $\beta = 1$ ($\COE$), $\beta=2$ ($\CUE$) and $\beta = 4$ ($\CSE$).
 \begin{defi}
For $\beta > 0$, let $\CbetN$ the distribution on $`T^N\times {\mathcal S}_N$ such that
%
\begin{itemize}
\item[1)] The joint law of $(\theta_1, \dots, \theta_N)$
 is :
\[\CbetN(\dd\!\theta_1, \cdots, \dd\!\theta_N) :=C_\beta(N)^{-1} |\Delta(\eul^{i\theta_1}, \cdots , \eul^{i\theta_n})|^\beta\ \dd\!\theta_1 \cdots \dd\! \theta_N\,,\]
where
$$C_\beta(N) = (2\pi)^N \frac{\Gamma\left(1 + \frac{\beta N}{2}\right)}{\left(\Gamma\left(1+\frac{\beta}{2}\right)\right)^N}.$$
\item[2)] The joint distribution of $(\w_1, \cdots , \w_N)$ on ${\mathcal S}_N$ is $\Dir_N(\beta/2)$.
\item[3)] The variables $\theta_{1},\cdots, \theta_{N}$ and $\w_{1},\cdots,\w_{N}$ are independent.
\end{itemize}
\end{defi}
Setting
\begin{equation}
\mu_\w\sn =\sum_{k=1}^N \w_k \delta_{\eul^{i\theta_k}}
\label{casse}
\end{equation}
\cite{Killip1} Killip and Nenciu proved in their Proposition 4.2 that under $\CbetN$ : 

1) the Verblunsky coefficients 
$c_j = c_j (\mu_\w\sn),(j=1,\ldots,N)$ are independent, 

2)  for $j=1,\ldots,N-1$,
\ben
\CbetN (c_j (\mu_\w\sn) \in \dd\! z) = \eta_{\frac{\beta}{2}(N-j) -1}(z) \dd\! z\,,
\een

3) the distribution of $c_{N}(\mu_\w\sn)$ is uniform on $\mathbb T$.

Actually, they built also  an explicit random matrix model (namely a five-diagonal matrix) 
whose spectral elements $(\eul^{i\theta_1}, \dots, \eul^{i\theta_N})$ and $(\w_1, \dots, \w_N)$ satisfy the previous properties. 
\subsection{$\so$ and Jacobi ensembles}
In this
subsection, we consider real random matrices. We first consider the subgroup $\so$ of $`U(2N)$, for which spectral measures (on $\mathbb T)$ are symmetric (i.e. invariant by complex conjugation). This entails that the  Verblunsky coefficients are real. We recall the results of \cite{Killip1} on the distributions of these coefficients, for $\so$  equipped with the Haar measure  and also for  a three-parameter family  of log-gases. 
Then, we project such spectral measures on $\mathbb R$  and,  by tuning the parameters we obtain random probability measures on $[0,1]$, whose moments are uniform. 
\subsubsection{$\so$ and extension}
If we provide $\so$ 
 with the normalized Haar measure $\lambda_{\so}$, then $\pm 1$ are a.s. not eigenvalues and  
  the spectral
measure 
 defined in (\ref{mastermu}) is  
  a.s. supported by pairwise conjugate complex numbers. It 
   may be written as:
\ben
\label{discsym}
\mu_\w^{(2N)} = 
\sum_{k=1}^N \frac{\w'_k}{2}\left(\delta_{\eul^{i\theta_k}} + \delta_{\eul^{-i\theta_k}}\right)\,.
\een
In that case, the distribution of $(\theta_1, \cdots, \theta_N)$ has a density proportional to
\ben
\label{deltatheta}
 |\Delta(2 \cos\theta_1, \cdots, 2\cos\theta_N)|^2 
\een
and the 
 array of weights $(\w'_1, \cdots , \w'_N)$ is
 independent
of $(\theta_1, \cdots, \theta_N)$ 
and satisfies
$$(\w'_1, \cdots , \w'_N) \el \Dir_N (1)\,,$$
(see \cite{Killip1} Proposition 3.4).\\
A natural generalization 
is obtained replacing the square in (\ref{deltatheta})
by an exponent $\beta$ as in Section \ref{goupil} and
adding a particular external potential:  
\begin{defi}
%
For $a,b, \beta > 0$ let $\wt J(\beta,a ,b, N)$ be the distribution of $(\theta_1, \cdots, \theta_N)$ having
a density proportional to
\ben
\label{convert}
|\Delta(\cos\theta_1, \cdots, \cos\theta_N)|^\beta \prod_{k=1}^N (1 - \cos \theta_k)^{a-\frac{1}{2}}(1 + \cos \theta_k)^{b-\frac{1}{2}}.
\een 
\end{defi}
In this  generalization, the 
random canonical moments 
have
some very remarkable properties  (like for the $\CbetN$ model) as soon as the weights are suitably distributed:
\begin{prop}(Killip-Nenciu \cite{Killip1} Prop. 5.3)
\label{53}
\begin{enumerate}
\item For $\mathbb{SO}(2N)$ equipped with the Haar measure $\lambda_{\so}$, the Verblunsky coefficients 
$c_1(\mu^{(2N)}_\w), \cdots , c_{2N-1}(\mu^{(2N)}_\w)$ are independent and satisfy:
\begin{eqnarray*}
c_{k}(\mu^{(2N)}_\w) &\el& \Beta_s\left(\frac{2N-k}{2}, \frac{2N-k}{2}\right)\ \ (k\leq 2N-1)\,.
\end{eqnarray*}
Moreover $c_{2N}(\mu^{(2N)}_\w)= -1$.
\item More generally, under the distribution $\wt J(\beta,a ,b, N)\otimes \Dir_N(\beta/2)$ on $[0, 2\pi]^N\times {\mathcal S}_N$, the Verblunsky coefficients $c_1(\mu^{(2N)}_\w), \cdots , c_{2N-1}(\mu^{(2N)}_\w)$ are independent and satisfy:
\begin{eqnarray*}
c_{k}(\mu^{(2N)}_\w) &\el& \Beta_s\left(\frac{2N-k-1}{4}\beta + a, \frac{2N-k-1}{4}\beta + b\right) \   k \ \hbox{odd} \  \\
c_{k}(\mu^{(2N)}_\w) &\el& \Beta_s\left(\frac{2N-k-2}{4}\beta + a +b , \frac{2N-k}{4}\beta \right)\   k \ \hbox{even} \,,
\end{eqnarray*}
for $k\leq 2N-1$. Moreover $c_{2N}(\mu^{(2N)}_\w)=-1$.
\end{enumerate}
\end{prop}
It is then natural to project the symmetric  measure $\mu_\w\sn$ on $\mathbb R$. It is the motivation of the next subsection.
\subsubsection{Jacobi ensembles} 
Let $\mu$ be a symmetric probability measure on $`T$ (i.e. invariant under complex conjugation) and let $\gamma = R(\mu)$  
be the pushforward of $\mu$ by the mapping $\eul^{i\theta} \mapsto \frac{1 +\cos \theta}{2}$, so that for $g$ continuous on $[0,1]$
\ben
\int_0^1 g(x) \dd\! \gamma (x) = \int_{`T} g\left(\frac{2+ z + \bar z}{4}\right)
\mu( \dd\! z)\,.
\een
Applying that to a probability measure with finite support
\[\mu_\w^{(N)} = \sum_{k=1}^N \frac{\w'_k}{2}\left(\delta_{\eul^{i\theta_k}} + \delta_{\eul^{-i\theta_k}}\right)\]
we get a probability measure
\ben
\label{34'}
\gamma_\w^{(N)} = R(\mu_\w^{(N)})= \sum_{k=1}^N \w'_k \delta_{x_k} \ \ \hbox{where} \  x_k = \frac{1 + \cos \theta_k}{2}\,, \ \ k=1, \dots, N\,.
\een
The correspondence between the Verblunsky/canonical coefficients of $\mu$ and the canonical coefficients of $\gamma = R(\mu)$ (in the sense of Subsection \ref{2.2}) is 
\ben
\label{palpha}                  
c_k (\gamma) = \frac{1}{2}\left(1 + c_{k}(\mu)\right) \ \left(\hbox{or, equivalently} \ c_k(\mu) = 2 c_{k}(\gamma) - 1\right)\,.
\een
This claim is proved by a geometrical argument in \cite{Simon2bis} (Theorem 13.3.1) and using Tchebychev polynomials in \cite{DS} (Section 9.3.8).
\begin{defi}
For $a,b,\beta>0$, let  $J(\beta,a ,b, N)$ the distribution on $[0,1]^{N}$  with density proportional to 
\[|\Delta(x_1, \dots, x_N)|^\beta\prod_{k=1}^N x_k^{b-1} (1-x_k)^{a-1}\,.\]
\end{defi}
If we equip $[0, 2\pi]^N$ with $\wt J(\beta,a ,b, N)$  
then the joint distribution of $x_{1},\cdots x_{n}$ is $J(\beta,a ,b, N)$.
We are now ready to claim :
\begin{lem}
\label{oneof}
Under the distribution $J(\beta,\frac{\beta}{4} ,\frac{\beta}{4}, N)\otimes \Dir_N(\beta/2)$,  the canonical moments  $c_1(\gamma_\w^{(N)}), \cdots , c_{2N-1}(\gamma_\w^{(N)})$ are independent and 
\ben
\label{loipk}
c_k (\gamma_\w^{(N)}) \el \Beta\left(\frac{2N-k}{4}\beta, \frac{2N-k}{4}\beta\right)\ \ (k \leq 2N-1)\,.
\een
\end{lem}
\proof From (\ref{palpha}) and (\ref{34'}) we have $c_k(\gamma_\w\sn) = \frac{1}{2}[1 + c_k(\mu_\w\sn)]$. Applying Proposition \ref{53} 2),
we see that the new coefficients inherit independence and have the mentioned distributions.\QED

\begin{theo}
\label{bizth}
\ 
\begin{enumerate}
\item If $\gamma_\w\sn = \sum_{k=1}^{N}\w_k \delta_{x_k}$ with
\ben
\label{biz1}
(x_1, \dots, x_{N},\w_1, \dots, \w_N)\el
J(4,1,1, N)\otimes \Dir_N(2) \,,\een
then the distribution of $\big(m_1 (\gamma_\w^{(N)}), \cdots, m_{2N-1}(\gamma_\w^{(N)})\big)$ is uniform on $M_{2N-1}^{[0,1]}$.
\item If $\gamma_\w\sn =  \w_0\delta_0 + \sum_{k=1}^{N}\w_k \delta_{x_k}$
\ben
\label{biz4}
(x_1, \dots, x_{N},\w_0, \dots, \w_N)\el  
J(4,1,3, N)\otimes \Dir (1, 2 , \cdots , 2, 2)\,,
\een
then the distribution of $\big(m_1 (\gamma_\w^{(N)}), \cdots, m_{2N}(\gamma_\w^{(N)})\big)$ is uniform on $M_{2N}^{[0,1]}$.
\item
If $\gamma_\w\sn = \w_0\delta_0 + \sum_{k=1}^{N-1}\w_k \delta_{x_k} + \w_N \delta_{1}$ with
\ben
\label{biz2}
(x_1, \dots, x_{N-1},\w_0, \dots, \w_N)\el
J(4,3,3, N-1)\otimes \Dir (1, 2 , \cdots , 2, 1)\,,\een
then the distribution of $\big(m_1 (\gamma_\w^{(N)}), \cdots, m_{2N-1}(\gamma_\w^{(N)})\big)$ is uniform on $M_{2N-1}^{[0,1]}$.
\item
If $\gamma_\w\sn =  \sum_{k=1}^{N}\w_k \delta_{x_k} + \w_{N+1} \delta_{1}$ with
\ben
\label{biz3}
(x_1, \dots, x_{N},\w_1, \dots, \w_{N+1})\el 
J(4,3,1, N)\otimes \Dir (2, 2 , \cdots , 2, 1)\,,\een
then the distribution of $\big(m_1 (\gamma_\w^{(N)}), \cdots, m_{2N}(\gamma_\w^{(N)})\big)$ is uniform on $M_{2N}^{[0,1]}$
\end{enumerate}
\end{theo}
To prove the theorem we  need the following proposition on principal representations of finite moment sequences (see for example \cite{Karlin1}, \cite{Karlin2} or \cite{DS} Defininition 1.2.10). These representations are related to the extreme values introduced in 
(\ref{doudou}) and (\ref{doudou1}).
\begin{prop}
\
\begin{enumerate}
\item Let $(m_{1},\cdots,m_{2N-1})\in M_{2N-1}^{[0,1]}$
\begin{itemize}
\item Lower principal representation: there exists a unique system $0<x_{1}<\ldots< x_{N}<1$ and $\w_{1},\ldots, \w_{N}\in(0,1)$  with $\sum_{k=1}^{N}\w_{k}=1$ such that the probability measure $\nu_{-}:=\sum_{k=1}^{N}\w_{k}\delta_{x_{k}}$ satisfies
$$m_{j}=\int_{0}^{1}x^{j}\nu_{-}(\dd\! x),\; j=1\ldots 2N-1.$$
\item Upper principal representation: there exists a unique system $0<x_{1}<\ldots< x_{N-1}<1$ and $\w_{0},\ldots, \w_{N}\in(0,1)$  with $\sum_{k=0}^{N}\w_{k}=1$ such that the probability measure $\nu_{+}:=\w_{0}\delta_{0}+\sum_{k=1}^{N-1}\w_{k}\delta_{x_{k}}+\w_{N}\delta_{1}$ satisfies
$$m_{j}=\int_{0}^{1}x^{j}\nu_{+}(\dd\! x),\; j=1\ldots 2N.$$
\end{itemize}
\item Let $(m_{1},\cdots,m_{2N})\in \Int(M_{2N}^{[0,1]})$
\begin{itemize}
\item Lower principal representation: there exists a unique system $0<x_{1}<\ldots< x_{N}<1$ and $\w_{0},\ldots, \w_{N}\in(0,1)$  with $\sum_{k=0}^{N}\w_{k}=1$ such that the probability measure $\nu_{-}:=\w_{0}\delta_{0}+\sum_{k=1}^{N}\w_{k}\delta_{x_{k}}$ satisfies
$$m_{j}=\int_{0}^{1}x^{j}\nu_{-}(\dd\! x),\; j=1\ldots 2N.$$
\item Upper principal representation: there exists a unique system $0<x_{1}<\ldots< x_{N}<1$ and $\w_{1},\ldots, \w_{N+1}\in(0,1)$  with $\sum_{k=1}^{N+1}\w_{k}=1$ such that the probability measure $\nu_{+}:=\sum_{k=1}^{N}\w_{k}\delta_{x_{k}}+\w_{N+1}\delta_{1}$ satisfies
$$m_{j}=\int_{0}^{1}x^{j}\nu_{+}(\dd\! x),\; j=1\ldots 2N.$$
\end{itemize}
\end{enumerate}
\label{thalys}
\end{prop}
\proof[of Theorem \ref{bizth}]
\begin{itemize}
\item To prove the first point,  take $\beta =4$ in Lemma \ref{oneof} and $n= 2N-1$ in Lemma \ref{unifreal}. 

\item To prove the remaining parts of the theorem we  use Proposition \ref{thalys}. We only treat the point 2) since the proof of the other ones may be tackled  identically. From Proposition  \ref{thalys} (even case, lower representation) the mapping between $(x_{1},\cdots,x_{N},\w_{1},\cdots, \w_{N})$ and $(m_{1},\cdots,m_{2N})$ is a diffeomorphism. Thus the probability density of the moments are completely known once we compute the Jacobian. Hence, we have to compute
$$
J =  \begin{vmatrix} x_1&\w_{1} & x_2& \w_{2}& \cdots &x_N& \w_{N}\\
x_1^2& 2\w_{1}x_1& x_2^2& 2\w_{2}x_2& \cdots &x_n^2& 2\w_{N}x_N\\
x_1^3& 3\w_{1}x_1^2& x_2^3& 3\w_{2}x_2^2& \cdots &x_N^3& 3\w_{N}x_N^2\\
\vdots& \vdots&\vdots& \vdots& &\vdots& \vdots\\ 
x_1^{2N}& 2N\w_{1}x_1^{2N-1}& x_2^{2N}& 2N\w_{2}x_2^{2N-1}& \cdots &x_{2N}^{2N}& 2N\w_{N}x_{2N}^{2N-1}  
\end{vmatrix}
$$
An easy calculation yields
$$
J = \left(\prod_{k=1}^{N}\w_k\right) \ \left(\prod_{k=1}^{N} x_k^2\right) \Delta(x_1, \cdots , x_{N})^4.
$$
Up to a normalizing constant, we recognize the density of $J(4,1,3, N)\otimes \Dir (1, 2 , \cdots , 2, 2)$. \QED
\end{itemize}
From Lemma \ref{unifreal}, the distribution of the canonical moments in the four cases of the previous theorem is 
completely known.
In their Section 2, Killip and Nenciu \cite{Killip1} give a model of tridiagonal real matrices admitting these spectral characteristics: their 
spectral random measures have the same distribution as 
$\gamma_\w$. Here, we present for the first and fourth cases of Theorem \ref{bizth} another interesting model issued from symmetric spaces. 
   
\begin{theo}
\label{unif2}Let $U$ be a $2n\times 2n$ random matrix such that $U \el g^Dg$ where $g$ is Haar distributed on $\mathbb S\mathbb O(2n)$. The spectral measure  $\mu_\w^{(n)}$ of $U$ is symmetric and  the distribution of  
$\big(m_1(\gamma_\w^{(n)}), \cdots , m_{n-1}(\gamma_\w^{(n)})\big)$   is
uniform on $M_{n-1}^{[0,1]}$.
\end{theo}

\proof Assume first that $n = 2N$. From 
Theorem \ref{bizth} 1.
  it is enough to check that this model induces the distribution $J(4,1,1,N)\otimes \Dir_N(2)$. We follow the notation of \cite{Duenez}. 

Let $\Phi$ be the mapping $g \in \mathbb S\mathbb O(2n) \mapsto g^D g$ and let ${\bf S}(n)$ be its image. Since  $\Phi^{-1}(1) := K(n) = \mathbb S\mathbb O(2n) \cap \USp(2n)$ we can see $\Phi$ as a one to one mapping from  $\mathbb S\mathbb O(2n)/ K(n)$ onto ${\bf S}(n)$. Actually $K(n)$ is isometric to $\mathbb U(n)$ via \[\mathbb U(n) \ni g \stackrel{i}{\mapsto} \begin{pmatrix}\Re g& -\Im g\\ \Im g&\Re g\end{pmatrix}\in K(n)\,. 
\]
Let $S$ be a random element of ${\bf S}(n)$ whose distribution is the pushforward of $\lambda_{\mathbb S\mathbb O(2n)}$ under $\Phi$.  We have
\[S \el i(g)a(i(g))^{-1}\]
where $g$ and $a$ are independent, 
$g$ is $\lambda_{\mathbb U (n)}$ distributed
 and if $n=2N$
\[a = \begin{pmatrix}\Re\Lambda_N& -\Im\Lambda_N&&\\
\Im\Lambda_N&\Re\Lambda_N&&\\
&&\Re\Lambda_N& \Im\Lambda_N\\
&&-\Im\Lambda_N&\Re\Lambda_N\end{pmatrix}\]
where $\Lambda_N = \diag (\eul^{i \theta_1}, \dots, \eul^{i\theta_N})$ and the $\theta_k$'s
have a joint density proportional to
\[\Delta(\cos\theta_1, \dots, \cos\theta_N)^4 \prod_{k=1}^N |\sin \theta_k|\,.\]
(see \cite{Duenez}, proof of 
 Theorem 2 (formulae 47 and 48)).
In terms of $x_k = (1+\cos \theta_k)/2$ this yields exactly the $J(4, 1,1, N)$ distribution. Moreover, 
a simple computation gives
\[\langle e_1, S^j e_1\rangle = \sum_{r=1}^N \w_r \cos (j\theta_r)\]
with $\w_r = (|g_{1,r}|^2 + |g_{1, N+r}|^2)$ and $(\w_1, \dots,\w_N)$ is $\Dir_N (2)$ distributed.

If $n = 2N+1$, similar considerations may be made to interpret formula (\ref{biz3}) as the random spectral measure of 
$\mathbb S\mathbb O(2n)/ K(n)$ when $n = 2N+1$ (in that case $1$ is double eigenvalue of $S$). 
\QED

\section{Large deviations for random spectral measures}
\label{arld}
We present here  the LDP for sequences of random measures defined in the previous section. The main result is
Theorem \ref{theogeneLoz} 
which is obtained as a consequence of our generic Proposition \ref{proponajim} concerning weighted random measures. This proposition states that the LDP for empirical measures with speed $N^2$ can be used to derive the LDP with speed $N$ for measures with 
 Dirichlet distributed weights.\\
For the sake  of completeness we briefly recall the LDP definition.
Let $(u_n)$
be a decreasing positive sequence of real numbers with $\lim_{n \rightarrow \infty}u_n=0$.
\begin{defi}
\label{dldp}
We say that a sequence $(R_{N})$ of probability measures on a measurable
Hausdorff space $(G,\br(G))$ satisfies the  LDP with rate function $I$ and speed
$(u_N^{-1})$ if:
\begin{itemize}
\item[i)] $I$ is lower semicontinuous (lsc), with values in
$\BBr^{+}\cup\{+\infty\}$.
\item[ii)] For any measurable set $A$ of $G$:
$$-I(\Int A)\leq
\liminf_{N\rightarrow\infty}u_N\log R_{N}(A)\leq
\limsup_{N\rightarrow\infty}u_N\log R_{N}(A)\leq
-I({\Clo A}),$$
where $I(A)=\inf_{\xi\in A}I(\xi)$ and $\Int A$ (resp. $\Clo A$) is the
interior (resp. the closure) of $A$.
\end{itemize}
We say that the rate function $I$ is good if its level set
$\{x\in G:\; I(x)\leq a\}$ is  compact for any $a\geq 0$.
More generally, a sequence of $G$-valued random variables 
is said to satisfy the
LDP if the sequence of their distributions satisfies the LDP.
\end{defi}
Hereafter, the space $G$ is the set of all non negative measures (or probability measures) supported by a
given  compact set. We endow this set with the topology of the weak convergence. 
The rate function obtained in this paper is be the so-called reversed Kullback information. 
 Let us recall that if $P$ and $Q$ be probability measures on $G$, the Kullback information between $P$ and $Q$ is defined by
\begin{equation}
{\mathcal K}(P | Q) = 
\begin{cases}
 \displaystyle\int_{G}\log\frac{\dd\! P}{\dd\! Q}\dd\! P\;\;\;\;&\mbox{if}\;P\ll Q\;\mbox{ and }\log\frac{\dd\! P}{\dd\! Q}\in L^1(P),\\
      +\infty\;\;&\mbox{ otherwise.}
\end{cases}
\end{equation}
For fixed $Q$, the nonnegative convex function ${\mathcal K}(\cdot | Q)$ is the rate function of the LDP for the empirical distribution of a sample according to $Q$
(Sanov's theorem \cite{DZ} p. 263).  
For fixed $P$, the function ${\mathcal K}(P | \cdot)$  is  nonnegative convex, it vanishes only for $Q=P$, and we call it reversed Kullback information with respect to $P$.  
\subsection{LDP for random measures with Dirichlet weights}
\label{varMEM}
\begin{prop}
\label{proponajim}
Let $\rho>1$ and $\{\xi_{k,N}, 1 \leq k\leq N\}_{N \geq 1}$ be a triangular array of random variables  belonging to a compact interval $[\chi_1,\chi_2]$ and let
$${\mathcal L}_N := \frac{1}{N} \sum_{k=1}^N \delta_{\xi_{k,N}}\,.$$
Assume that the sequence $({\mathcal L}_N)$ 
satisfies the LDP in ${\mathcal M}_1([\chi_1, \chi_2])$ with speed $(N^\rho)$ with a good rate function $I_\xi$. 
Assume further that $I_\xi$ has a unique minimum at $\nu$ whose support  is 
$[\chi_1,\chi_2]$.
Let  $\{\w_{k,N}, 1 \leq k\leq N\}_{N \geq 1}$ another triangular array, independent of the first one, such that for $N \geq 1$, the vector $(\w_{k,N}, 1 \leq k\leq N)$ follows the distribution $\Dir_N (a)$ ($a>0$).  Then, the family
\ben
\mu_{N} :=  \sum_{k=1}^N \w_{k,N} \delta_{\xi_{k,N}}
\een
 satisfies the LDP  in ${\mathcal M}_1([\chi_1, \chi_2])$ with speed $(N)$ and  good rate function 
is $a{\mathcal K}(\nu | \cdot)$.
\end{prop}

\proof

\noindent 
The proof  uses heavily the classical representation (\ref{10}) :
\ben
\label{repres}
\{\w_{k,N}, 1 \leq k\leq N\} \el \left(\frac{Y_1}{Y_1 + \cdots +Y_N}, \cdots, \frac{Y_N}{Y_1 + \cdots +Y_N}\right)
\een
where the $Y$'s are independent and $\gamma_a$ distributed. The modified measure
\[\wt\mu_N := \frac{1}{N}\sum_{k=1}^N Y_k \delta_{\xi_{k,N}}\,,\]
has then independent weights and is easier to handle. We come back to the original measure with
\[\mu_\w\sn = \frac{\wt\mu_N}{\wt\mu_N(1)}\,.\]
Now, recall
that the cumulant generating function of the $\gamma_a$ distribution is
$$
L_Y(\tau): =
\begin{cases}
a \log(1-\tau)\;\;\;\;&\mbox{if}\; \tau<1,\\
    +\infty\;\;&\mbox{ otherwise,}
\end{cases}
$$
and that its 
Cram\'er transform 
is (see for example \cite{DZ}) 
$$
j_a(x) = 
\begin{cases}
 x-a -a \log \frac{x}{a}\;\;\;\;&\mbox{if}\; x>0,\\
    +\infty\;\;&\mbox{ otherwise.}
\end{cases}
$$
To understand the flavor of our proof, let us first assume that the sequence $(\xi_{k,N})$ 
is not random, and that the sequence $({\mathcal L}_N)$ converges weakly to $\nu$.
In that case, using the main Theorem of
\cite{Najim1} 
we see that the sequence $\wt\mu_N$ satisfies the LDP with speed $(N)$ and good rate function :
\ben
\wt I (\mu):= \int_{[\chi_1,\chi_2]}(j_a(g)-g)\dd\! \nu + \int_{[\chi_1,\chi_2]}\dd\! \mu\;\; (\mu\in{\mathcal M}([\chi_1,\chi_2]))\,,
\een
where
$\dd\!\mu=g \dd\!\nu+(\mu-g \dd\!\nu)$ is the Lebesgue decomposition of $\mu$ with respect to $\nu$.
Further,  as almost surely $\wt\mu_N([\chi_1,\chi_2])>0$ and using the contraction principle
(see \cite{DZ} p.126), we may conclude that $(\mu_N)$ satisfies the LDP with speed $(N)$ and good rate function  
\begin{eqnarray*}
I (\xi) = \inf \{\wt I (\mu):\;\mu\in{\mathcal M}([\chi_1,\chi_2]),\;  \mu = \mu([\chi_1,\chi_2]) \xi \} = \inf_{b > 0} \wt I(b\xi).
\end{eqnarray*}
Therefore, a direct computation yields 
\[I(\xi) = a {\mathcal K}(\nu | \xi)\,.\]
Let us now turn to the general case.
We have to get rid of the randomness of the $\xi_{k,N}$. But, we  show
that, since the LDP with larger speed holds for its empirical measure, this
randomness has no effect on the rate function. To get 
 the LDP,
we 
use a method inspired by the proof of the so-called Gärtner-Ellis-Baldi theorem 
(\cite{DZ} p.157). Roughly speaking, it consists of two steps.  
First we compute the limiting normalized generating function and its convex conjugate function. Then we conclude by density of exposed points of this last function.

To begin, 
 let $f$ be a continuous function on $\BBr$ such that $L_Y\circ f$ is bounded. Then, integrating first on the random variables $Y$ we may write
\ben
\nonumber
I_N := {\mathbb E} \exp\wt\mu_{N}(f) = 
  {\mathbb E}\exp\Big[\sum_k L_Y(f(\xi_{k,N}))\Big] = {\mathbb E}\exp \big(N {\mathcal L}_N(L_Y\!\circ\!f)\big)\,. 
\een
Now, from the assumption, the sequence of random variables $({\mathcal L}_N(L_Y\!\circ\!f))$
converges to $\int (L_Y\!\circ\!f)\!\ \dd\! \nu$ and satisfies the LDP with speed
$(N^\rho)$. 
Now fix $\varepsilon > 0$ and set
$$A := \{\eta \in {\mathcal M}_+ ([\chi_1, \chi_2]) : \left|\eta(L_Y\!\circ\!f)  - \nu(L_Y\!\circ\!f) \right|\leq \varepsilon\},$$
and $I_N = I_{N, A} + I_{N, A^c}$ with 
\ben
I_{N, A}:= {\mathbb E}\left[ \exp\big( N  {\mathcal L}_N(L_Y\!\circ\!f)\big) \mathbf{1}_{\{{\mathcal L}_N \in A\}}\right].
\een
Obviously, we have
\ben
\label{estim1}
\hspace{1cm}
\exp[N\left(\nu (L_Y\!\circ\!f)  - \varepsilon \right)] {\mathbb P}\left({\mathcal L}_N \in A\right) \leq I_{N, A}\leq \exp [N\left(\nu (L_Y\!\circ\!f)  + \varepsilon \right)]
\een
and
\ben
\label{estim2}
I_{N, A^c}\leq \exp [N \Vert L_Y\!\circ\!f\Vert_\infty] \ {\mathbb P}\left({\mathcal L}_N \in A^c\right)
\een
The LDP upper bound gives 
\ben
\label{2maj}
\limsup \frac{1}{N^\rho} \log {\mathbb P}({\mathcal L}_N \in A^c) \leq -\inf \{I_\xi (\eta) ; \eta \in \overline {A^c}\} 
\,.
\een
Since $I_\xi$ is good, its infimum on the closed set $\overline {A^c}$ is reached and since it has a unique global minimum at $\nu$ which is not in $\overline {A^c}$, the right side of (\ref{2maj}) is negative. Since $\rho > 1$ we get
$
\limsup_N \frac{1}{N} \log I_{N, A^c} = -\infty,
$
and then
\ben
\label{majo}
\limsup_N \frac{1}{N} \log I_{N} \leq \nu(L_Y\!\circ\!f)  + \varepsilon\,.
\een
To get the lower bound, we have
\begin{eqnarray*}
\liminf_N \frac{1}{N} \log I_{N} \geq \liminf_N \frac{1}{N} \log I_{N, A}\geq\nu (L_Y\!\circ\!f)  - \varepsilon + \liminf_N \frac {1}{N} \log {\mathbb P}\left({\mathcal L}_N \in A\right)
\end{eqnarray*}
From (\ref{2maj}) it is clear that  $\lim {\mathbb P}\left({\mathcal L}_N \in A\right) = 1$ so that
\ben
\label{mino}
\liminf_N \frac{1}{N} \log I_{N} \geq \nu(L_Y\!\circ\!f)  - \varepsilon
\een
As (\ref{majo}) and (\ref{mino}) hold for every $\varepsilon > 0$, we may conclude
$$\lim_N \frac{1}{N} \log I_{N} = \nu(L_Y\!\circ\!f) \,.$$
The remaining of the proof is the same as the LDP proofs developed in \cite{grz} (see also
\cite{gamb}).
 \QED

\subsection{Large deviations for  spectral measures in the circular and Jacobi ensembles}
We now state the LDP for the sequence of random probability measures $(\mu_\w\sn)$ 
in the above frameworks. Notice that 
in the case $\beta = 2$, with the help of Theorem \ref{=C}, the first result has yet been shown  
in \cite{LozEJP} and the second one in \cite{FabLoz} Theorem 2.3. In this last paper, the proofs use the so-called projective limit approach. We  take here a completely different way. Indeed, we give a general proof by using a variation on the results developed in \cite{Najim1} for random measures. As a matter of fact, we have shown in Subsection \ref{varMEM} the LDP for random measures
 with weights having a Dirichlet distribution.  
\begin{theo}
\label{theogeneLoz}
\
\begin{enumerate}
\item Assume that for any $N$, $\mathbb{T}^N\times\mathcal{S}_{N}$ is endowed with the $\CbetN$ distribution. Recall that $\mu_\w\sn$ is built  as in (\ref{casse}).
Then, the sequence of random probability distributions $(\mu_\w\sn)$  
satisfies the LDP with speed $N$  
and  good rate function
\ben
I(\xi)= \frac{\beta}{2}{\mathcal K}\left(\lambda_{`T} | \xi \right),\;
(\xi\in{\mathcal M}_1(`T))
\een
where $\lambda_{`T}$ is the uniform distribution on $`T$. 
\item  Let $a, b > 0$ be fixed. The sequence  $(\xi^{(N)})$ under the $J(\beta, a, b , N)\otimes \Dir_N(\beta/2)$ distribution satisfies the LDP with speed $N$ and good rate function
\ben
I (\xi) = \frac{\beta}{2}{\mathcal K}(\ARCSIN | \xi)
\een
\end{enumerate}
\end{theo}
\proof
The proof of this theorem is a direct consequence of Proposition \ref{proponajim}. Indeed, from \cite{HiaiIHP} it is known that the empirical distribution built on $(\theta_k)$ 
satisfies the LDP with speed $N^2$ whose rate function has a unique minimum in $\lambda_{\mathbb T}$.
%
For the Jacobi ensemble, Hiai and Petz (\cite{hiai2} Section 2) proved the LDP with speed $N^2$ 
for the empirical measure built on $(x_k)$. Moreover, they have shown that - with our assumptions -  
the unique minimizer of the rate function for this LDP is
the  $\ARCSIN$ distribution.
\begin{rem}
The same conclusion holds true under the probability measures of Theorem \ref{bizth}.
\end{rem}
\begin{rem} The assumptions about the supports are crucial. 
 In a companion paper \cite{gr+} we study important models in which they are not satisfied.  For the $\GUE_N$, 
 the support of ${\mathcal L}_N$ is not included in a fixed compact set and the limiting distribution is the semicircle. 
In this  example, we  find an additional term in the rate function. 
For the  $J(\beta, \tau_1 N, \tau_2 N ; N)$, 
 the support of ${\mathcal L}_N$ is included in $[0,1]$ but the limiting distribution  is supported by a compact subinterval of $(0,1)$. It is called the Kesten-MacKay distribution and generalizes the $\ARCSIN$ distribution. 
In \cite{gr+}, we prove the LDP but the explicit form of the rate function is only conjectured. The same is true for the Laguerre ensemble.
\end{rem}
\begin{rem}
Other discrete random measures with Dirichlet weights appear in the literature. 
Let $\alpha$ 
be  a positive measure on a compact space $K$. A Dirichlet process $\mu$ of parameter $\alpha$ is a 
probability distribution denoted by ${\mathcal D}(\alpha)$ on ${\mathcal M}_1(K)$. It is  such that for any measurable finite partition of $K$, $(A_1, \cdots, A_N)$, we have  
$
(\mu(A_1) , \cdots, \mu(A_N)) \el \Dir_N (\alpha(A_1) , \cdots, \alpha(A_N)).
$\\
If $\alpha$ is the finite discrete measure $\theta \sum_{k=1}^N \delta_{y_k}$
where 
 $y_1,\cdots, y_N\in K,$ we have 
\ben
\label{mesdir2}
{\mathcal D}\left(\theta \sum_{j=1}^N \delta_{y_j}\right) \el \sum_{k=1}^N Y_k \delta_{y_k}
\een
where 
$
(Y_1, \cdots, Y_N) \el \Dir_N (\theta)\,. 
$
This means that the random measures studied in the previous sections are mixed Dirichlet processes. The mixing distribution is the law of the random eigenvalues. 
Ganesh and O'Connell (\cite{Ganesh2}) study  
  ${\mathcal D}(\alpha + \sum_{i=1}^N \delta_{y_i})$ when 
 $N^{-1}\sum_{i=1}^N \delta_{y_i}$ converges to some measure $\nu$ and when $\alpha$ is a measure whose
support is $K$. They show that $\mu$
satisfies the LDP with speed $N$, and good rate function
$
I(\cdot) = {\mathcal K}(\nu | \cdot)$.
In  \cite{DawsonFengII} there is an extension with an infinite number of random locations, 
(see also \cite{DawsonFeng2} for the connection with  the Poisson-Dirichlet distribution).
\end{rem}
\subsection{Relation with spherical integrals}
In the unitary model, the direct computation of the limiting cumulant generating functional of  $\mu_\w\sn$ leads to a spherical integral.  Indeed, if $\varphi\in C(\mathbb T)$ (the set of all continuous function on $\mathbb T$),
we have \[\mu_\w\sn(\varphi) = \langle e_1, \varphi(U)e_1\rangle = \tr(G_NVD\sn_\varphi V^*)
\,\] 
where  
$D\sn_\varphi := \diag (\varphi(\eul^{i\theta_1}), \cdots, \varphi(\eul^{i\theta_N}))\,,$ (the $\eul^{i\theta_k}$ are the eigenvalues of $U$), $V \in \mathbb U(N)$ 
 and $G_N := \diag(1, 0, \cdots, 0)$. 
Consequently, the Laplace transform of $\mu_\w$ is, for $\varphi\in C(\mathbb T)$, 
\ben\mathbb E \exp\big(N\mu_\w\sn (\varphi)\big)  =
\mathbb E_\theta\!\ I^{(2)}_N(G_N , D_\varphi^{(N)})
%
\label{cheval}
\een 
where
$$I^{(2)}_N(G_N , D_\varphi^{(N)}) = \left(\int_{\mathbb U(N)} \exp [N \tr(G_NVD\sn_\varphi V^*)] \ \dd\! V\ \right),$$
($\mathbb E_\theta$ denotes here expectation with respect to the variables $\theta_1, \dots, \theta_N$).
It is interesting to notice that this last spherical integral over the unitary group can be expressed as a hypergeometric function with two matrix arguments (see for example \cite{Orlov} p. 97).

In \cite{AliceM} and \cite{CollinsS} it is proved that  
$
N^{-1} \log I_N^{(2)}(G_N, D_\varphi\sn) 
$
has a limit as $N \rightarrow \infty$, as soon as  the empirical spectral distribution of $D_\varphi^{(N)}$
converges weakly.
More precisely  in \cite{AliceM} Th.6 (see also \cite{CollinsS} Th. 4) it is proved that
\ben
\lim_N \frac{1}{N} \log I_N^{(2)}(G_N, D_\varphi\sn) = F_\varphi(1),
\een
where $F_\varphi(1)$  
is defined in the following way. Let $\nu_\varphi$ be the limit of 
 the empirical spectral distribution of $D\sn_\varphi$, (in our case it is the image by $\varphi$ of $\Lambda_{`T}$). 
Further let $\varphi_{\min}:=\min_{z\in `T}\varphi(z)$ and 
 $\varphi_{\max}:=\max_{z\in `T}\varphi(z)$.
The Stieltjes transform of $\nu_\varphi$ is defined for $x\in (-\infty,\varphi_{\min})\cup (\varphi_{\max},+\infty)$ by:
\ben
H_\varphi (x) = \int_\BBt \frac{1}{x- \varphi(z)} \lambda_{_\BBt}(\dd\! z)\,.
\een
We have that $H^\downarrow_\varphi := \lim_{x\downarrow \varphi_{\max}} H_\varphi(x) < 0 < H^\uparrow_\varphi := \lim_{x\uparrow \varphi_{\min}} H_\varphi(x)$  and the range of $H_{\varphi}$ is
$(H^\downarrow_\varphi,0)\cup (0, H^\uparrow_\varphi)$.
 (\cite{AliceM} Property 9).
The  Voiculescu $R$-transform is the function $R_\varphi$ satisfying, for all $y$ in the range of $H_{\varphi}$, 
\ben
H_\varphi\left(R_\varphi(y) + \frac{1}{y}\right) = y\,.
\een
As a result (\cite{AliceM} Theorem 6),
\ben
\nonumber
F_\varphi(1) &=&v(1) - \int_\BBt \log (1 + v(1) -\varphi(z)) \  \lambda_{_\BBt}(\dd\! z)\\
  v(1)&=&
     \begin{cases}
        R_\varphi(1) & \text{if $H_\uparrow \leq 1 \leq H_\downarrow$} \\
        \varphi_{\max} - \displaystyle 1 & \text{if $1 > H_\downarrow$}\\
        \varphi_{\min} - \displaystyle 1 & \text{if $1 < H_\uparrow$}
     \end{cases}
\een
Now, the 
LDP for the empirical distribution built on $\theta_{1},\cdots, \theta_{N}$ holds with speed $N^2$. So that,  $N^{-1} \ln \mathbb E \exp \big(N\mu_\w\sn(\varphi)\big)$ 
and $N^{-1} \log I_N^{(2)}(G_N, D_\varphi\sn)$ have the same limit (to
show this claim just use the continuity of  spherical integrals proved in \cite{Maida} Prop. 2.1).
The rate function of the LDP for $\mu_\w\sn$ can be recovered by 
taking the supremum in $\varphi \in C(\mathbb T)$ of 
$\mu(\varphi) -F_\varphi(1)$, where $\mu \in {\mathcal M}_1(\mathbb T)$. 
Setting $g(\eul^{i\theta}) = \varphi(\eul^{i\theta}) - v(1)$, 
 we have
\ben
\label{GMT}
\mu(\varphi) - F_{\varphi}(1) = 
\int_{\mathbb T} g(z)  \mu(\dd\! z)  + \int_{\mathbb T}\log(1-g(z))\  \lambda_{\mathbb T}( \dd\! z)\,.
\een
Taking the supremum in 
$g\in C(\mathbb T)$, we recover the well-known duality formula
\[ \sup_{g \in C(\mathbb T)}\left[ \int_{\mathbb T} g(z) d\mu(z)  + \int_{\mathbb T}\log(1-g(z))\  \lambda_{\mathbb T}(\dd\! z)\right] = {\mathcal K}(\lambda_{\mathbb T} | \mu)\,.\]
References and several consequences of this formula for the associated moment problem may be found in \cite{FabLoz}.
%
\\
\ \\
\noindent {\bf Acknowledgments.}
The authors are very grateful to the 
anonymous referees for a careful
reading of the paper. A.R. was partially funded by the Agence Nationale de la
Recherche grant ANR-08-BLAN-0311-03.

\bibliographystyle{plain}

\bigskip
$$\;$$
Fabrice Gamboa\\
Universit\'e de Toulouse\\ 
Universit\'e Paul Sabatier\\
Institut de Math\'ematiques de Toulouse\\
F-31062 Toulouse, France\\
{\tt fabrice.gamboa@math.univ-toulouse.fr}
$$\;$$
Alain Rouault         
LMV  B\^atiment Fermat\\
 Universit\'e
Versailles-Saint-Quentin \\
F-78035 Versailles\\ France\\
{\tt  Alain.Rouault@math.uvsq.fr}
\end{document}